\documentclass[a4paper,12pt]{amsart}
\usepackage{amssymb}
\numberwithin{equation}{section}
\newtheorem{theorem}{Theorem}[section]
\newtheorem{prop}[theorem]{Proposition}
\newtheorem{lem}[theorem]{Lemma}
\newtheorem{cor}[theorem]{Corollary}
\newtheorem{defi}[theorem]{Definition}

\theoremstyle{remark}
\newtheorem{rem}[theorem]{Remark}
\newcommand{\R}{\mathbb{R}}

\newcommand{\N}{\mathbb{N}}

\newcommand{\ds}{\displaystyle}

\renewcommand{\bar}{\overline}

\author[J.~Benameur and Ch. Elhechmi  ]{Jamel Benameur and Chokri Elhechmi}

\address{Research Laboratory Mathematics and Applications LR17ES11, Department of Mathematics, Faculty of Science of Gab\`es, university of Gab\`es; Tunisia}
\email{\sl jamelbenameur@gmail.com}
\email{\sl chokri.elhechmi@issatm.rnu.tn}

\title[Iterative method for solving a nonlinear Robin problem]
{Iterative method for solving a nonlinear Robin problem}

\begin{document}
	\begin{abstract}
		In \cite{CJ1} M. Jaoua et al. studied the linear approximation of Robin problem on $\Omega$ an open bounded domain of $\R^d$, and they given some important results. In this paper, we study a nonlinear approximation of an elliptic problem with a nonlinear Robin boundary condition in a domain of $\R^2$. We prove the existence and uniqueness of solution by an iterative construction method with admissible condition on $\partial\Omega$.
	\end{abstract}

	%@@@@@@@@@@@@@@@@@@@@@@@@@@@@@@@@@@@@%@@@@@@@@@@@@@@@@@@@@@@@@@@@@@@@@@@@@%@@@@@@@@@@@@@@
	

	%@@@@@@@@@@@@@@@@@@@@@@@@@@@@@@@@@@@@%@@@@@@@@@@@@@@@@@@@@@@@@@@@@@@@@@@@@%@@@@@@@@@@@@@@
	\maketitle
	\tableofcontents

	%@@@@@@@@@@@@@@@@@@@@@@@@@@@@@@@@@@@@%@@@@@@@@@@@@@@@@@@@@@@@@@@@@@@@@@@@@%@@@@@@@@@@@@@@@
	
\section{\bf Introduction}\label{sec-intro}

In this paper we are interested in the study of an elliptic P.D.E. defined on an open  bounded domain  of $\R^2$ with a nonlinear Robin type boundary condition defined on a part of the boundary $\Gamma_R \subset \partial \Omega$ where $\partial \Omega$ is  the boundary of $\Omega$. It is well known that P.D.E. models large number of scientific fields such as  physical, chemical, biological phenomenon,.... For example, referring to M. Vogelius et al. \cite{vxu}, this equation can provide information  about the corrosion effect that can attack a metal plate. Indeed, the authors considered a metallic plate, occupying the $\Omega\subset\R^d$ domain,\,$d\geq2$, attacked by corrosion on part of its boundary. The electrostatic potential satisfies the potential equation $\Delta u = 0 $ in $\Omega$ and the corrosion effect is modeled by the following equation:
 \begin{equation}
 \displaystyle \frac{\partial u}{\partial n} = \lambda (e^{\alpha u} - e^{-(1-\alpha )u}),
 \end{equation}
where $$\frac{\partial u}{\partial n} = \nabla u . n :\;{\rm denotes\, the\, normal\, derivative,}$$
$\lambda$ is an exchange coefficient and $\alpha \in]0, 1[$ is a kinetic parameter.\\
The Taylor expansion of order 2 of this condition gives us
\begin{equation}
 \displaystyle \frac{\partial u}{\partial n} \approx \lambda (u + \displaystyle\frac{2 \alpha -1 }{2} u^2).
 \end{equation}
This problem was the subject of several research works. J. M. Arrieta and S. M. Bruschi  \cite{JMAr-SMBr}   study how oscillations in the boundary of a domain affect the behavior of solutions of elliptic equations with nonlinear boundary conditions of the type $$\displaystyle\frac{\partial u}{\partial n}  +g(x,u) = 0$$ when the boundary of the domain presents a highly oscillatory behavior.\\

 In  his paper  \cite{Pasc-Cher}  P. Cherier consider a nonlinear Neumann problems on $C^{\infty}$ compact riemannian manifolds $(\overline{\Omega}, g)$ of dimension $n \geq 2$ whose boundary $\partial \Omega$ is an $(n - 1)$ dimensional sub-manifold. He study a nonlinear Robin problem of the form : $$ \Delta u + f(u,x) = 0 \;{\rm in}\;\Omega$$ and
 $$\displaystyle\frac{\partial u}{\partial n}+g(u,y) = 0\;{\rm on}\;\partial \Omega.$$
 Under the assumptions   the functions $|f(t,x)|$ and $|g(t,y)|$ are bounded by $C(1 +|t|^\sigma)$ or
 $C \exp(|t|^\sigma)$ for all $(t,x,y) \in \R \times \Omega \times \partial \Omega$, the author proved the existence of a regular  solution $u \in C^{\infty}(\overline{\Omega})$.\\

The linear case, when  the parameter $\alpha = \frac{1}{2}$, was the subject of several studies. For example, M. Jaoua et al (see \cite{CJ1}) are studied an inverse problem of identifying the coefficient $\lambda$ by measurements on the boundary. Our aim in this paper is the study of existence and  uniqueness of the following elliptic problem where the Robin condition is non linear. Also we construct an iterative solution of this problem:
$$
(RP)\left\{\begin{array}{l}
\displaystyle\Delta u= 0\;{\rm in}\;\Omega\\
\displaystyle u= 0\;{\rm on}\;\Gamma_D\\
\displaystyle\frac{\partial u}{\partial n}= \phi(x)\;{\rm on }\;\Gamma_N\\
\displaystyle\frac{\partial u}{\partial n}+\varphi(x)\big(u+\frac{2\alpha-1}{2}u^2)= g(x)\;{\rm on}\;\Gamma_R,
\end{array}\right.$$
 where $\Gamma_D$, $\Gamma_N$ and $\Gamma_R$ are {\bf three open subsets} of $\partial\Omega$ such that
$$(H1)\left\{\begin{array}{l}
\partial\Omega=\overline{\Gamma_D}\cup \overline{\Gamma_N}\cup\overline{\Gamma_R}\\
\Gamma_D\cap \Gamma_N=\Gamma_D\cap\Gamma_R=\Gamma_N\cap\Gamma_R=\emptyset\\
\sigma(\Gamma_D)>0,\;\sigma(\Gamma_N)>0,\;\sigma(\Gamma_R)>0,
\end{array}\right.
$$
with $\sigma $ is the  superficial measure on $\partial\Omega$  induced  by  the  Lebesgue  measure on $\overline{\Omega}$.\\
The functions  $\varphi,\;\psi$ are two continuous functions from $\overline{\Gamma_R} $ to $\R$. The function  $\phi$ (the current flow imposed on the part of the boundary $\Gamma_N$)   verifies $\phi \in L^2(\Gamma_N) $ and $ \phi \not \equiv 0$.\\

The paper is organized as follows. In Section \ref{sec-linPr}, we present the linear problem and we recall some results relating to this problem.  In section \ref{sec-prlem}, we collect some preliminaries results and we give some notation. In section \ref{sec-GnPrb}, we formulate the general nonlinear problem and we prove our main results: existence, uniqueness and the iterative construction of the solution of the nonlinear problem. The Section  \ref{sec-Appendix} is an appendix composed of some definitions and main properties that we use in this work.
\section{\bf The linear Robin problem} \label{sec-linPr}
Before treating the non-linear case, we give a quick reminder about the linear case. Precisely, we consider the following system
$$
(P_0)\left\{\begin{array}{l}
\displaystyle\Delta u= 0\;{\rm in}\;\Omega\\
\displaystyle u= 0\;{\rm on}\;\Gamma_D\\
\displaystyle\frac{\partial u}{\partial n}= \phi(x)\;{\rm on } \;\Gamma_N\\
\displaystyle\frac{\partial u}{\partial n}+\varphi(x)u= g(x)\;{\rm on}\;\Gamma_R.
\end{array}\right.
$$
To write the result of \cite{CJ1}, we need some definitions and notations, let
$$
V =\{ v \in H^1(\Omega);\;  {\rm such\,  that }\;  v = 0\;{\rm on }\; \Gamma_D \}.
$$
 $V$  is a Hilbert space for the scalar product $  \langle u, v\rangle_V  = \ds \int_\Omega  \nabla u. \nabla v$.\\
The associate norm is defined by
$$\| u\|_{V} =\Big(\int_{\Omega} |\nabla u|^2\Big)^{1/2} =\|\nabla u\|_{L^2(\Omega)}\, \forall u \in V.$$
Let  he trace map  defined by
  $$
 \begin{array}{llll}
 \tau_1 &: V &\rightarrow &H^{1/2}(\Gamma_R)\\
  &u &\mapsto & u_{|\Gamma_R}.
 \end{array}
 $$
The map $\tau_1$ is continuous linear operator if the space $H^{1/2}(\Gamma_R)$ is equipped by the norm of $L^2(\Gamma_R)$, then  there exists a positive constant $C$ such that :
 \begin{equation}\label{eq1}
 \|u_{|\Gamma_R}\|_{L^2(\Gamma_R)}   \leq C \|u\|_{V} \,\,\, \forall u \in V.
\end{equation}
 We denote $$\beta_1:=\inf\big\{ C>0;\;\|u_{|\Gamma_R}\|_{L^2(\Gamma_R)}\leq C \|u\|_{V} \,\,\, \forall u \in V\big\}=\||\tau_1\||.$$
The  problem $(P_0)$ is equivalent to the following variational formulation:
$$\left\{\begin{array}{l}
{\rm Find }\;u\in V\;{\rm such\; that}\\\\
a(u,v)=l(v),\;\forall v\in V,
\end{array}\right.$$
where
$$\begin{array}{lcl}
a(u,v)&=&\displaystyle\int_{\Omega}\nabla u \nabla v  +\int_{\Gamma_R}\varphi uv\\
l(v)&=&\displaystyle\int_{\Gamma_N}\phi v+\int_{\Gamma_R}g v.
\end{array}$$
The set of admissible coefficients $\Phi_{ad}$ (associate to problem $(P_0)$) is defined by
  $$\Phi_{ad} =\{ \varphi \in C^0(\bar{\Gamma_R});\;{\rm such\, that }\; \displaystyle\inf_{x\in\bar{\Gamma_R}} \varphi(x) > -\displaystyle\frac{1}{\beta_1^2}\}.$$
 
\begin{theorem}\label{th01}(\cite{CJ1})  Let $\varphi\in\Phi_{ad}$. Then the problem  $(P_0)$ has a unique solution $u_0\in V$ and
\begin{equation}\label{eq2}
\|u_0\|_{V}  \leq  \ds\frac{\beta_1}{1+C_\varphi\beta_1^2}\big(\|\phi\|_{L^2(\Gamma_N)}  + \|g\|_{L^2(\Gamma_R) }\big),
\end{equation}
where \begin{equation}\label{eqC1}C_\varphi = \min (\ds\inf_{\Gamma_R} \varphi, 0)\leq0.\end{equation} Moreover, if $1+C_\varphi\beta_1^2\geq1/2$, we get
\begin{equation}\label{eq3}
\|u_0\|_{V}  \leq  \ds 2\beta_1\big(\|\phi\|_{L^2(\Gamma_N)}+ \|g\|_{L^2(\Gamma_R) }\big).
\end{equation}
\end{theorem}

\begin{rem} To prove this theorem, the authors use a classical methods. They start by proving that the bilinear function $a$ is continuous, coercive , and the linear function $l$ is continuous in the space $V$. Then, by applying the Lax Milgram theorem they deduce the existence  and uniqueness of a solution to the problem  $(P_0)$.  The inequality (\ref{eq3}) is obtained by using the trace map $\tau_1$.

\end{rem}
%\begin{theorem}\label{th02}(\cite{CJ1}) Let $\varphi_1,\varphi_2\in\Phi_{ad}$ and $u_1,\;u_2$ be the solutions of problem $(P_0)$ with $\varphi_1$ (resp. $\varphi_2$) as a Robin coefficient. Suppose that ${u_1}_{|K}={u_2}_{|K},$ then $\varphi_1=\varphi_2$.
%\end{theorem}
\vspace{0.5cm}
Our aim in this paper is the study of existence and  uniqueness of the elliptic problem $(RP)$, where the Robin condition is non linear. Also we construct an iterative solution of this problem. To study problem $(RP)$ and to find the principle result, we need some notations, definitions and preliminary results
 which will be presented in the following paragraph.
\section{\bf Preliminaries results and main result}\label{sec-prlem}
In this section, we collect some classical definitions, and we give
some lemmas which are well suited to the study of the Robin problem $(RP)$.
The proof of this result is based on Lax Milgram's theorem, then it is necessary to prove that  the considered bi-linear function is coercive and continuous.
Also we will prove that the considered linear function is continuous in the Hilbert space $V$.

Let $p\in [1, \infty[$, for a bounded and regular domain $\Omega$,  we define the following Sobolev spaces :
$$W^{1,p}(\Omega)=\{u\in L^p(\Omega);\;\nabla u\in L^p(\Omega)\}.$$
%-------
\begin{lem}\label{lem01}(\cite{DJ-Mer})  Let $\Omega\subset\R^2$ be a Lipschitzian domain and $p\in(1,\infty)$,then
%$$W^{1,p}(\Omega)\hookrightarrow W^{1-\frac{1}{p},p}(\partial\Omega).$$
$$\|u_{|\partial\Omega}\|_{W^{1-\frac{1}{p},p}(\partial\Omega)}\leq C_{p,\Omega} \|u\|_{W^{1,p}(\Omega)}$$
where the  constant $C_{p,\Omega} > 0$ depend only on $p$ and domain $\Omega$.
\end{lem}
\begin{lem}\label{lem02}(\cite{Jer-Dron}) Let $\Omega\subset\R^2$ be a Lipschitzian domain. Then, for $q\in[2,+\infty)$ we have
$$H^{1/2}(\partial\Omega)\hookrightarrow L^q(\partial\Omega).$$
Precisely, there is a constant
$C_{\Omega,q} > 0$, depending on $q$ and domain $\Omega$, such that
$$\|u_{|\partial\Omega}\|_{L^q(\partial\Omega)}\leq C_{\Omega,q} \|u\|_{H^{1/2}(\partial\Omega)}.$$
\end{lem}
\begin{lem}\label{lem03}(H\"older inequality)
Let $(X, \mathcal A, \mu)$ be a measure space.
\begin{enumerate}
\item Let $p, q \in [1,\infty]$ with $1/p + 1/q = 1$. Then for all measurable real- or complex-valued functions $f$ and $g$ on $X$ such that $f \in L^p(X)$ and $g \in L^q(X)$, we have $fg\in L^1(X)$ and
$$\|fg\|_{L^1(X)}\leq \|f\|_{L^p(X)}\|g\|_{L^q(X)}.$$
\item Let $p$, $q$, $r \in [1, \infty[$ satisfy $\ds \frac{1}{p} + \frac{1}{q} = \frac{1}{r}$. For all
    $f \in L^p(X)$ and $g \in L^q(X)$ we have $fg \in   L^r(X)$ and
$$
\ds \|fg\|_{L^r(X)} \leq \|f\|_{L^p(X)}\|g\|_{L^q(X)}.
$$
\end{enumerate}
\end{lem}
\begin{cor}\label{cor01} Let $1<p<\infty$. For all $u, v \in L^p(X)$, we have  $u |v|^{p-1}  \in   L^1(X)  $ and we have
$$\int_X\Big| u(x) |v(x) |^{p-1}\Big|d\mu(x)\leq \|u\|_{L^p(X)}\|v\|_{L^p(X)}^{p-1}.$$
%Particularly, if $p=k\in\N,\;p\geq2$, we have
%$$|\int_Eu(x)v(x)^{k-1}dx|\leq \|u\|_{L^k(E)}\|v\|_{L^k(E)}^{k-1}.$$
\end{cor}
\begin{rem} In the same way of definition of $\tau_1$, we define the trace map $\tau_2$ by
 $$
 \begin{array}{llll}
 \tau_2 &: V &\rightarrow &H^{1/2}(\Gamma_N) \\
  &u &\mapsto & u_{|\Gamma_N}
 \end{array}
 $$
 and by using the lemma \ref{lem02} we have
$$
 \|u_{|\Gamma_N}\|_{L^3(\Gamma_N)} \leq C \|u\|_{V} \,\,\, \forall u \in V.
$$
 We denote
 \begin{equation}\label{eq11}
 \beta_2=\inf\{ C>0;\;\|u_{|\Gamma_N}\|_{L^3(\Gamma_N)}\leq C \|u\|_{V} \,\,\, \forall u \in V\}.
\end{equation}
 \end{rem}
\noindent In order to state our main result we need some notations and hypothesis on the boundary conditions:
$$
\begin{array}{lll}
&{\bf (C1)}& \phi\in L^2(\Gamma_N)\;{\rm and}\;g\in L^2(\Gamma_R)\\
&{\bf (C2)}& M_0=2\beta_1(\|\phi\|_{L^2(\Gamma_N)}+\|g\|_{L^2(\Gamma_R)})\\
&{\bf (C3)}& {\rm For}\;\varepsilon\;{\rm be\; real\; number\; satisfying}: 0<\varepsilon<\inf\displaystyle\Big(\frac{1}{4\beta_1^2},\frac{1}{2\beta_2^3M_0}\Big),\;{\rm we\; notice}\\
&&\Phi_{ad}^\varepsilon=\{\varphi\in C(\overline{\Gamma_R});\;\|\varphi\|_{L^\infty(\Gamma_R)}<\varepsilon\}.
\end{array}$$
Precisely, we have the following result.
\begin{theorem}\label{thJC01} Let $\varphi\in\Phi_{ad}^\varepsilon$.
\begin{enumerate}
\item The problem $(RP)$ has a solution $u$ in $V$, moreover $u$ satisfying
\begin{equation}\label{eq3-10}
\|u\|_{V}  \leq M_0.
\end{equation}
\item This solution is unique in $\overline{B_V(0,M_0)}=\{h\in V:\;\|h\|_V\leq M_0\}.$
\end{enumerate}
\end{theorem}
\begin{rem}${}$
\begin{enumerate}
\item The proof of this theorem will be in the following section where we will study a more general case and this theorem will be a particular one.
\item We can prove the uniqueness in $\overline{B_V(0,R)}=\{h\in V:\;\|h\|_V\leq R\}$ for each $R\geq M_0$, but we must take some smallness conditions on $\varepsilon$, which imposes some smallness on $\|\varphi\|_{L^\infty(\Gamma_R)}$.
\end{enumerate}
\end{rem}
\section{\bf Proof of the main result}\label{sec-GnPrb}
In this section we study  a general Robin problem. Precisely, we give a complete study of the following nonlinear Robin problem :
$$
(GRP)\left\{\begin{array}{l}
\displaystyle\Delta u= 0\;{\rm in}\;\Omega\\
\displaystyle u= 0\;{\rm on}\;\Gamma_D\\
\displaystyle\frac{\partial u}{\partial n}= \phi(x)\;{\rm on }\;\Gamma_N\\
\displaystyle\frac{\partial u}{\partial n}+\varphi(x)u+\psi(x)u^2=g(x)\;{\rm on}\;\Gamma_R.
\end{array}\right.
$$
Our approach consists of proving the existence and uniqueness of the problem $(GRP)$ by using an iterative method allowing the construction of the unique solution of this problem. Before starting this study, it is important to define  a set of admissible functions:
$$\begin{array}{lll}
&{\bf(C4)}&\;{\rm For}\;\varepsilon_1,\varepsilon_2>0\;{\rm such\; that}\;\varepsilon_1<\displaystyle\frac{1}{4\beta_1^2}\;{\rm and}\; \varepsilon_2<\displaystyle\frac{1}{4\beta_2^3M_0},\;{\rm we\; define}\\
&&\Phi_{ad}^{\varepsilon_1,\varepsilon_2}=\bigg\{(\varphi,\psi)\in C(\overline{\Gamma_R})^2;\;\inf_{\Gamma_R}\varphi>-\varepsilon_1\;{\rm and}\;\|\psi\|_{L^\infty(\Gamma_R)}<\varepsilon_2\bigg\}.\end{array}$$
We now give the main result in this paragraph :
\begin{theorem}\label{thJC21}
Let $(\varphi,\psi)\in\Phi_{ad}^{\varepsilon_1,\varepsilon_2}$.
\begin{enumerate}
\item The problem $(GRP)$ has a solution $u$ in $V$ satisfying the following inequality
\begin{equation}\label{eq3-1}
\|u\|_{V}  \leq M_0.
\end{equation}
\item This solution is unique in $\overline{B_V(0,M_0)}=\big\{h\in V:\;\|h\|_V\leq M_0\big\}.$
\end{enumerate}
\end{theorem}
\begin{rem}${}$
\begin{enumerate}
\item The proof of the existence result is done by classical method:
\begin{itemize}
\item Approximate system $(P_k)$.
\item Existence and uniqueness of solution $u_k$ of the approximate system $(P_k)$ with some uniform estimates.
\item Convergence of the sequence $(u_k)_{k\geq0}$ to a function $u\in B_{M_0}$.
\item Passage to the limit in the approximate system, which implies that $u$ is a solution of $(GRP)$.
\item The uniqueness is inspired by the proof of the existence step.\\
\end{itemize}
\item We can prove the uniqueness in $B_R=\big\{h\in V:\;\|h\|_V\leq R\big\}$ for each $R\geq M_0$, but we must take some smallness conditions on $\varepsilon_1$ and $\varepsilon_2$.
\end{enumerate}
\end{rem}
\subsection{Some preliminaries Lemmas}
 In order to prove this theorem, we need the following lemmas which we will prove firstly.
\begin{lem}\label{lem11} Let $f\in V$ such that $\|f\|_V\leq M_0$. Then the bilinear application
$$A_f:V\times V\rightarrow\R,\;(u,v)\mapsto \int_{\Omega}\nabla u.\nabla v+\int_{\Gamma_R} \varphi u
+\displaystyle\int_{\Gamma_R} \psi fuv  $$
is continuous and coercive. Precisely, we have
$$\begin{array}{lcl}
|A_f(u,v)|&\leq&C_2\|u\|_V\|v\|_V,\;\forall (u,v)\in V^2\\\\
A_f(u,u)|&\geq&\displaystyle\frac{1}{2}\|u\|_V^2,\;\forall u\in V,
\end{array}$$
where $C_2 = 1+\beta_1^2\|\varphi\|_{L^{\infty}(\Gamma_R)}  +\beta_2^3M_0 \|\psi\|_{L^{\infty}(\Gamma_R)}$.
\end{lem}
{\bf Proof.} We demonstrate the continuity and coercivity of the application $A_f$ separately.
\begin{enumerate}
\item[$\bullet$] {\bf Continuity of $A_f$:} For $(u,v) \in V \times V$ and by using H\"older inequality, we get
$$
\begin{array}{lcl}
 A_f(u,v) &=&\displaystyle\int_{\Omega}\nabla u(x).\nabla v(x)dx+\int_{\Gamma_R} \varphi u v +\int_{\Gamma_R} \psi fu v \\\\
&\leq&\displaystyle\| \nabla u\| _{L^2(\Omega)}\| \nabla v\|_{L^2(\Omega)} + \|\varphi\|_{L^{\infty}(\Gamma_R)}  \|u\|_{L^2(\Gamma_R)}\|v\|_{L^2(\Gamma_R)}\\
&&+\displaystyle\|\psi\|_{L^{\infty}(\Gamma_R)} \|f\|_{L^3(\Gamma_R)}\|u\|_{L^3(\Gamma_R)}\|v\|_{L^3(\Gamma_R)}\\\\
&\leq&\displaystyle\|u\| _{V}\|v\|_{V} + \|\varphi\|_{L^{\infty}(\Gamma_R)}  \beta_1^2\|u\|_{V}\|v\|_{V}\\
&&+ \|\psi\|_{L^{\infty}(\Gamma_R)} \alpha_2^3\|f\|_{V}\|u\|_{V}\|v\|_{V}\\\\
&\leq& \big(1+\beta_1^2\|\varphi\|_{L^{\infty}(\Gamma_R)}  +\beta_2^3M_0 \|\psi\|_{L^{\infty}(\Gamma_R)}\big) \|u\|_{V}\|v\|_{V}\\\\
&\leq& C_2 \|u\|_{V}\|v\|_{V}.
\end{array}$$
\item[$\bullet$] {\bf Coercivity of $A_f$:} For $u\in V$ and by using H\"older inequality \ref{lem03}, we get
$$\begin{array}{lcl}
 A_f(u,u) &=&\ds\| \nabla u\|^2_{L^2(\Omega)} +\displaystyle\int_{\Gamma_R} \varphi u^2 +\displaystyle\int_{\Gamma_R} \psi fu^2\\
&\geq& \|u\|^2_V+\displaystyle\inf_{\Gamma_R}\varphi \int_{\Gamma_R}u^2  - \|\psi\|_{L^{\infty}(\Gamma_R)}\int_{\Gamma_R} |f|u^2\vspace{3mm}\\
&\geq&  \| u\|^2_V + \ds \inf_{\Gamma_R}  \varphi \|u\|^2_{L^2(\Gamma_R)}  - \|\psi\|_{L^{\infty}(\Gamma_R)} \|f\|_{L^3(\Gamma_R)}\|u\|^2_{L^3(\Gamma_R)}.
\end{array}$$
By using the equations (\ref{eq1}),(\ref{eqC1}) and (\ref{eq11}), we obtain $$
A_f(u,u)\geq\displaystyle\|u\|^2_V+ \beta_1^2\displaystyle C_\varphi \|u\|^2_V-\beta_2^3 \|\psi\|_{L^{\infty}(\Gamma_R)} M_0  \|u\|^2_V.$$
By using the fact that $(\varphi,\psi)\in \Phi_{ad}^{\varepsilon_1,\varepsilon_2}$, we get
$$\begin{array}{lcl}
A_f(u,u)&\geq&\displaystyle\|u\|^2_V-\varepsilon_1\beta_1^2\|u\|^2_V-\beta_2^3 \varepsilon_2 M_0  \|u\|^2_V\\
&\geq&\displaystyle\bigg(1  - \varepsilon_1\beta_1^2-\beta_2^3 \varepsilon_2 M_0  \bigg)\|u\|^2_V.
\end{array}$$
The choices of $\varepsilon_1$ and $\varepsilon_2$ (See {\bf(C4)}) imply
\begin{eqnarray*}
A_f(u,u) \geq \frac{1}{2} \|u\|^2_{V},\;\forall u \in V.
\end{eqnarray*}
\end{enumerate}
\begin{lem}\label{lem22}
Let $\phi \in L^2(\Gamma_N)$ and $g \in L^2(\Gamma_R)$ then the linear function
$$
l : V\rightarrow\R,\;v \mapsto \int_{\Gamma_N}\phi v +\int_{\Gamma_R} g v
$$
is continuous in the Hilbert space $V$.
\end{lem}
{\bf Proof:} For $v\in V$, we have
$$
\begin{array}{lcl}
|l(v)|& =&\displaystyle|\int_{\Gamma_N}\phi v+\int_{\Gamma_R}g v|\\
&\leq&\displaystyle|\int_{\Gamma_N}\phi v|+|\int_{\Gamma_R}g v|.
\end{array}
$$
By applying the Cauchy-Schwarz inequality and using the continuity of $\tau_1$, we get
$$
\begin{array}{lcl}
|l(v)|& \leq&\displaystyle\| \phi\|_{L^2(\Gamma_N)} \|v\|_{L^2(\Gamma_N)}+  \|g\|_{L^2(\Gamma_R)}\|v\|_{L^2(\Gamma_R)}\\
 &\leq&\displaystyle\beta_1\big(\| \phi\|_{L^2(\Gamma_N)}+ \|g\|_{L^2(\Gamma_R)}\big)\|\nabla v\|_{L^2(\Omega)}\\
 &\leq& C_0\beta_1\| v\|_{V},
\end{array}
$$
where $C_0  = \| \phi\|_{L^2(\Gamma_N)}+ \|g\|_{L^2(\Gamma_R)}$.
\begin{rem}
It is clear that  if $\varphi\in \Phi_{ad}^{\varepsilon}$, then $ (\varphi,(\alpha-\frac{1}{2})\varphi))\in \Phi_{ad}^{\varepsilon_1,\varepsilon_2}$, with  $$\varepsilon_1=\varepsilon<\frac{1}{4\beta_1^2},\;\varepsilon_1=\frac{\varepsilon}{2}<\frac{1}{4\beta_2^3M_0}.$$
    Then Theorem \ref{thJC21} implies Theorem \ref{thJC01}.
\end{rem}
%--------
\subsection{Proof of Theorem \ref{thJC21}}
\subsubsection{\bf Existence} The proof is done in three steps:
\begin{enumerate}
\item[$\bullet$]{\bf First step:} In this section we give a general result of intermediate linear Robin system.
\begin{prop}\label{thprop1} Let $\phi\in L^2(\Gamma_D)$, $g\in L^2(\Gamma_R)$, $(\varphi,\psi)\in\Phi_{ad}^{\varepsilon_1,\varepsilon_2}$ and $f\in V$ such that
$$\|f\|_V\leq M_0.$$
Then, the following problem
$$
({\rm IR}(f))\left\{\begin{array}{l}
\displaystyle\Delta u= 0\;{\rm in}\;\Omega\\
\displaystyle u= 0\;{\rm on}\;\Gamma_D\\
\displaystyle\frac{\partial u}{\partial n}\;= \phi(x)\;{\rm on }\;\Gamma_N\\
\displaystyle\frac{\partial u}{\partial n}+\varphi(x)u+\psi(x)f(x)u\;=g(x)\;{\rm on}\;\Gamma_R
\end{array}\right.
$$
has a unique solution $u$ in $V$. Moreover, we have
\begin{equation}\label{eq31}
\|u\|_{V}  \leq  M_0.
\end{equation}
\end{prop}
{\bf Proof.} The variational formulation of the problem $({\rm IR}(f))$ is given by
$$
A_f(u,v) = l(v),\,\, \forall \, (u, v) \in V\times V.
$$
where
$$\begin{array}{lcl}
A_f(u,v)&=&\displaystyle\int_{\Omega}\nabla u(x)\nabla v(x)dx +\int_{\Gamma_R}\varphi u v +\int_{\Gamma_R}\psi fuv \\
l(v)&=&\displaystyle\int_{\Gamma_N}\phi v+\int_{\Gamma_R}g v.
\end{array}$$
By Lemmas \ref{lem11}-\ref{lem22}, we have $A_f$ is continuous, coercive and $l$ is continuous. Then by applying Lax-Milgram Theorem and taking into account that $(\varphi,\psi)$ is in $\Phi_{ad}^{\varepsilon_1,\varepsilon_2}$, there is a unique $u\in V$ solution of $({\rm IR}(f))$. Moreover, we have
$$
\|u\|_V\leq 2\beta_1(\|g\|_{L^2(\Gamma_R)}+\|\phi\|_{L^2(\Gamma_N)})=M_0.
$$
\item[$\bullet$]{\bf Second step:} In this step, we give an approximate schema to problem $(GRP)$. Let $u_0\in V$ the unique solution of $(P_0)$ given by Theorem \ref{th01}. Particularly, $u_0$ satisfies $\|u_0\|_V\leq M_0$.\\
For each   $k\in\N$, we consider the following problem :
$$
(P_k)\left\{\begin{array}{l}
\displaystyle\Delta u= 0\;{\rm in}\;\Omega\\
\displaystyle u= 0\;{\rm on}\;\Gamma_D\\
\displaystyle\frac{\partial u}{\partial n}= \phi(x)\;{\rm on }\;\Gamma_N\\
\displaystyle\frac{\partial u}{\partial n}+\varphi(x)u+\psi(x)u_{k-1}u=g(x)\;{\rm on}\;\Gamma_R.
\end{array}\right.
$$
For $k=1$ and by applying the first step(Proposition \ref{thprop1}) with $f=u_0$, we deduce that the problem  $(P_1)=({\rm IR}(u_0))$ has a unique solution $u_1$ in $V$.\\
 Moreover, we have
\begin{equation}\label{eq31step2}
\|u_1\|_{V}  \leq2\beta_1\big(\|\phi\|_{L^2(\Gamma_N)}  + \|g\|_{L^2(\Gamma_R) }\big)=M_0.
\end{equation}
For $k\in\N, \;k\geq2$, suppose that, we construct $u_1,...,u_k\in V$ such that for each $1\leq i\leq k-1$, $u_{i+1}\in V$ is the unique solution of the problem
$$
(P_{i+1})\left\{\begin{array}{l}
\displaystyle\Delta u= 0\;{\rm in}\;\Omega\\
\displaystyle u= 0\;{\rm on}\;\Gamma_D\\
\displaystyle\frac{\partial u}{\partial n}=\phi(x)\;{\rm on }\;\Gamma_N\\
\displaystyle\frac{\partial u}{\partial n}+\varphi(x)u+\psi(x)u_{i}u=g(x)\;{\rm on}\;\Gamma_R.
\end{array}\right.
$$
and suppose that
\begin{equation}\label{eq31step2i}
\|u_{i+1}\|_{V}  \leq M_0.
\end{equation}
Again, by the first step we can construct $u_{k+1}$ the unique solution of $(P_{k+1})=({\rm IR}(u_k))$ satisfying $\|u_{k+1}\|_{V}  \leq M_0.$\\
Finally, we construct a sequence $(u_k)_{k\in\N}$ of element of $V$ such that $u_k$ is the unique solution of $(P_k)$ and $\|u_k\|\leq M_0$.
\item[$\bullet$]{\bf Third step:}  In this step, we will prove that  the sequence $(u_k)_{k\in\N}$ converges in $V$ to a solution of $(GRP)$. Precisely, we have the following result.
\begin{prop}\label{thprop3} There is $u\in V$ such that $\|u\|_V\leq M_0$ and
\begin{equation}\label{eq33}
\lim_{k\rightarrow+\infty}\|u_k-u\|_{V}=0.
\end{equation}
Moreover, $u$ is a solution of $(GRP)$ in $V$.
\end{prop}
{\bf Proof.} Let $w_k= u_{k+1}- u_k$, where $u_k$ and $u_{k+1}$ are respectively the solutions of the problem $(P_k)$ and $(P_{k+1})$. Then, we have
$$\left\{\begin{array}{l}
\displaystyle\Delta w_{k}= 0 \;{\rm in}\;\Omega\\
\displaystyle w_{k}= 0\;{\rm on}\;\Gamma_D\\
\displaystyle\frac{\partial w_{k}}{\partial n}= 0\;{\rm on }\;\Gamma_N\\
\displaystyle\frac{\partial w_{k}}{\partial n}+\varphi(x)w_{k}+\psi(x)u_kw_{k}=-\psi(x)u_k w_{k-1}\;{\rm on}\;\Gamma_R.
\end{array}\right.
$$
\noindent$\bullet$ {\bf We begin by proving the following lemma.}
\begin{lem} There exists a positive constant $\mathcal{K}< 1$ such that
\begin{equation}\label{eq14}
\forall \, k \in \N,\;\;\|  w_k \|_{V}  \leq \mathcal{K} \| w_{k-1} \|_{V}.
\end{equation}
\end{lem}
{\bf Proof of lemma :}  Let $v\in V$. By using the Green formula, we obtain
$$
\begin{array}{lcl}
\displaystyle\int_{\Omega}  \nabla w_k \nabla v &=& -\ds \int_{\partial \Omega} \frac{\partial w_k}{\partial n} v\\
&=&\displaystyle-\int_{\Gamma_R} \varphi w_k v - \int_{\Gamma_R }\psi u_k w_k v - \ds\int_{\Gamma_R }\psi u_k w_{k-1} v.
\end{array}$$
Then
$$
\displaystyle\int_{\Omega}  \nabla w_k \nabla v +\int_{\Gamma_R} \varphi w_k v +\int_{\Gamma_R }\psi u_k w_k v = -\displaystyle\int_{\Gamma_R} \psi u_k w_{k-1} v.
$$
If we take $v=w_k$, we get
$$
\begin{array}{lcl}
I_k &:=&\displaystyle \int_{\Omega} |\nabla w_k|^2 + \int_{\Gamma_R }\varphi |w_k|^2 + \int_{\Gamma_R} \psi u_k |w_k|^2\\
&=&- \displaystyle\int_{\Gamma_R} \psi u_k w_{k-1} w_k\\
&\leq&\displaystyle\|\psi\|_{\infty} \|u_k\|_{L^3(\Gamma_R)}\|w_{k-1}\|_{L^3(\Gamma_R)}\|w_{k}\|_{L^3(\Gamma_R)}\\
&\leq&\displaystyle\beta_2^3\|\psi\|_{\infty} \|u_k\|_{V}\|w_{k-1}\|_{V}\|w_{k}\|_{V}\\
&\leq&\displaystyle M_0\beta_2^3\varepsilon_2\|w_{k-1}\|_{V}\|w_{k}\|_{V}.
\end{array}
  $$
In other hand, by using the continuity of the trace functions $\tau_1$ and $\tau_2$ and H\"older inequality, we obtain
$$
\begin{array}{lcl}
I_k &=&\displaystyle \int_{\Omega} |\nabla w_k|^2 + \int_{\Gamma_R} \varphi |w_k|^2 + \int_{\Gamma_R} \psi u_k |w_k|^2\\
&\geq&\displaystyle\|w_{k}\|_{V}^2+(\inf_{\Gamma_R}\varphi)\|w_k\|_{L^2(\Gamma_R)}^2-\int_{\Gamma_R}|\psi|. |u_k| |w_k|^2\\
&\geq&\displaystyle\|w_{k}\|_{V}^2+(\inf_{\Gamma_R}\varphi)\|w_k\|_{L^2(\Gamma_R)}^2-\|\psi\|_{L^\infty(\Gamma_R)}\| u_k\|_{L^3(\Gamma_R)}\|w_{k}\|^2_{L^3(\Gamma_R)}\\
 &\geq&\displaystyle\|w_{k}\|_{V}^2 -\varepsilon_1 \beta_1^2 \|w_{k}\|^2_{V} - \varepsilon_2 \| u_k\|_{L^3(\Gamma_R)}\|w_{k}\|^2_{L^3(\Gamma_R)}\\
 &\geq&\displaystyle\|w_{k}\|_{V}^2 -\varepsilon_1 \beta_1^2 \|\nabla w_{k}\|^2_{L^2(\Omega)} - -\varepsilon_2\beta_2^3\| u_k\|_{V}\|w_{k}\|^2_{V}\\
 &\geq&\displaystyle\|w_{k}\|_{V}^2 -\varepsilon_1\beta_1^2 \|w_{k}\|^2_{V} - \varepsilon_2 M_0 \beta_2^3\|w_{k}\|^2_{V}\\
   &\geq&\displaystyle\big(  1 -\varepsilon_1\beta_1^2 - \varepsilon_2 M_0 \beta_2^3\big)\|w_{k}\|^2_{V}\\
    &\geq&\displaystyle\frac{1}{2}\|w_{k}\|^2_{V}.
\end{array}
$$
 Hence
 \begin{eqnarray} \label{eq4}
 \|w_{k}\|_{V} \leq \mathcal{K} \|w_{k-1}\|_{V}.
 \end{eqnarray}
 where
 \begin{eqnarray} \label{eq5}
 \mathcal{K} =\displaystyle2\beta_1^2\varepsilon_2 M_0.
 \end{eqnarray}
By applying the conditions of {\bf(C4)}, we  deduce that  $\mathcal{K} < 1$.\\\\

{\bf Let's return to the proof of Theorem \ref{thJC21}},  by combining this result with the fact that $V$ is an Hilbert space, we get: the series $\sum_{k\geq1}w_k$ converges in $V$ to an element $S\in V$. Therefore, the sequence $(u_k)$ solutions of $(P_k)$ converges in the space $V$ to an element $u\in V$, precisely
$$\displaystyle u=\lim_{k\rightarrow\infty} u_{k+1} = S+u_0  \mbox{  in }  V.$$\\
\noindent$\bullet$ {\bf Prove that $u$ is a solution of $(GRP)$.}\\
We have $\displaystyle u=\lim_{k\rightarrow+\infty} u_{k}$ in $V$ and $u_k$ is a solution of $(P_k)$ for each $k$, then
$$
\begin{array}{l}
\displaystyle\Delta u=0\;{\rm in}\;\Omega \\
\displaystyle u=0\;{\rm on}\;\Gamma_D\\
\displaystyle\frac{\partial u}{\partial n}=\phi(x)\;{\rm on}\;\Gamma_N.
\end{array}
$$
It remains to show that
$$\frac{\partial u}{\partial n}+\varphi(x)u+\psi(x)u^2=g(x)\;{\rm on}\;\;\Gamma_R.$$
For the linear parts, we have
$$\lim_{k\rightarrow+\infty}\frac{\partial u_k}{\partial n}+\varphi(x)u_k=\frac{\partial u}{\partial n}+\varphi(x)u\;{\rm on}\;\;\Gamma_R.$$
For the nonlinear part, we have
$$\psi(x)u^2-\psi(x)u_{k-1}u_k=\psi(x)u(u-u_k)+\psi(x)u_k(u_{k-1}-u)\;{\rm on}\;\;\Gamma_R.$$
By applying Cauchy-Schwarz and H\"older inequalities and the continuity of the trace functions $\tau_1$ and $\tau_2$, we get
$$\begin{array}{lcl}
\displaystyle\int_{\Gamma_R}|\psi(x)u^2-\psi(x)u_{k-1}u_k|
&\leq&\displaystyle\int_{\Gamma_R}|\psi(x)|.|u|.|u-u_k|+\int_{\Gamma_R}|\psi(x)|.|u_k|.|u_{k-1}-u|\\
&\leq&\displaystyle\|\psi\|_{L^\infty(\Gamma_R)}\Big(\int_{\Gamma_R}|u|.|u-u_k|+\int_{\Gamma_R}|u_k|.|u_{k-1}-u|\Big)\\
&\leq&\displaystyle\varepsilon_2\Big(\int_{\Gamma_R}|u|.|u-u_k|+\int_{\Gamma_R}|u_k|.|u_{k-1}-u|\Big)\\
&\leq&\displaystyle\varepsilon_2\Big(\|u\|_{L^2(\Gamma_R)}\|u-u_k\|_{L^2(\Gamma_R)}
+\|u_k\|_{L^2(\Gamma_R)}\|u_{k-1}-u\|_{L^2(\Gamma_R)}\Big)\\
&\leq&\displaystyle\varepsilon_2\Big(\beta_1^2\|u\|_{V}\|u-u_k\|_{V}
+\beta_1^2\|u_k\|_{V}\|u_{k-1}-u\|_{V}\Big).
\end{array}$$
By using the fact that
$$\|u\|_{V}\leq M_0\;{\rm and}\;\|u_k\|_{V}\leq M_0,\;\forall k\geq1,$$
we get
$$\int_{\Gamma_R}|\psi(x)u^2-\psi(x)u_{k-1}u_k|\leq \varepsilon_2\beta_1^2M_0\Big(\|u-u_k\|_{V}
+\|u_{k-1}-u\|_{V}\Big).$$
By using the fact that $\displaystyle u=\lim_{k\rightarrow+\infty} u_{k}$ in $V$, we get
$$\ds\lim_{k\rightarrow+\infty}\int_{\Gamma_R}|\psi(x)u^2-\psi(x)u_{k-1}u_k|=0.$$
Which complete the proof and $u$ is a solution of $(GRP)$ in $V$.
\end{enumerate}
\subsubsection{\bf Uniqueness} Let $u,v\in \overline{B_V(0,M_0)}$ two  solutions of the problem $(GRP)$ such that $u$ is given by the last step, and let us prove that $v=u$.\\
For this, put $w=u-v$ and by applying the same calculus for $w_k$ to $w$, we get
$$\|w\|_V\leq \mathcal K\|w\|_V.$$
Then $(1-\mathcal K)\|w\|_V\leq0$. As $0<\mathcal K<1$, we get $\|w\|_V=0$ and $w=0$. \\
Therefore $v=u$, and the uniqueness in $\overline{B_V(0,M_0)}$ is proved.
\section{\bf Appendix} \label{sec-Appendix}
\subsection{Appendix A: The Hilbert space $V$}
\begin{lem}
The space $V$ is a Hilbert space for the scalar product $\langle u, v\rangle_V  = \ds \int_\Omega  \nabla u. \nabla v$ where
$$
V =\{ v \in H^1(\Omega);\;  {\rm such\,  that }\;  v = 0\;{\rm on }\; \Gamma_D \}.
$$
\end{lem}
{\bf Proof.} Let $(u_n)$ a Cauchy sequence in $V$. Since the norm  $\|.\|_V$ and the norm $|.|_{H^1(\Omega)}$ are equivalent then the sequence $(u_n)$ is a Cauchy sequence in the  Hilbert space $H^1(\Omega)$.\\
Let  $u^*$ be the limit of $(u_n)$ in $H^1(\Omega)$. \\
By using the continuity of trace map $\tau : H^1(\Omega) \rightarrow H^{1/2}(\partial\Omega)$ we deduce
$$
{u_n}_{|_{ \partial \Omega}} \longrightarrow u^*_{|_{ \partial \Omega}}\;{\rm in}\;H^{1/2}(\partial\Omega)
$$
hence $u^* = 0$ in $\Gamma_D$.\\
From where  $u^*\in V$. As a result $V$ is a complete space.
\subsection{Appendix B: Sobolev spaces  $W^{1,p}(\Omega)$}
\begin{defi}
Let $\Omega \subset \R^{d}$ be an open domain,  $1\leq p\leq\infty$ and $m\in\N^{*}$, we define the following
Sobolev space
\begin{equation}
 W^{m,p}(\Omega)=\big\{u\in L^{p}(\Omega),  D^{\alpha}u\in L^{p}(\Omega), \forall \alpha\in\mathbb{N}^d,  |\alpha|\leq m\big\},
 \end{equation}
with the associate norm
$$\|u\|_{W^{m,p}(\Omega)}\, =\, \|u\|_{L^{p}(\Omega)}\,+\,\sum^{m}_{\alpha=1}\|D^{\alpha}u\|_{L^{p}(\Omega)}.$$
\end{defi}
\begin{rem}
\item If $p=2$ we define $H^{m}(\Omega)\, = \,W^{m,2}(\Omega)$ by
\begin{equation}
H^{m}(\Omega)=\big\{u\in L^{2}(\Omega),\, D^{\alpha}u\in L^{2}(\Omega),\, \alpha\in\mathbb{N}^d,\, |\alpha|\leq m\big\}.
\end{equation}
\item $H^{m}(\Omega)$ is a Hilbert space for the scalar product
$$(u,v)_{H^{m}(\Omega)}\,=\,(u,v)_{L^{2}(\Omega)}\,+\,\sum^{m}_{|\alpha|=1}(D^{\alpha}u,D^{\alpha}v),$$
with the norm
  $$\|u\|_{H^{m}(\Omega)}\, =\Big(\|u\|_{L^{2}(\Omega)}^2\,+\,\sum^{m}_{|\alpha|=1}\|D^{\alpha}u\|_{L^{2}(\Omega)}^2\Big)^{1/2}.$$
\end{rem}
\begin{defi} For $p \in [2, \infty[$ we define the trace space
$$
W^{1-\frac{1}{p},p}(\partial\Omega)= \tau\big(W^{1,p}(\Omega)\big)  = \{ \tau(u), u \in W^{1,p}(\Omega)\}
$$
with the norm
$$
\|f\|_{W^{1-\frac{1}{p},p}(\partial\Omega)} = \inf \{ \|u\|_{W^{1,p}(\Omega)} \, \mbox{ such that } \tau(u) =f\},
$$
where $\tau$ is the map trace : $\tau : W^{1,p}(\Omega) \rightarrow L^{p}(\partial\Omega),\;u\mapsto u_{|\partial\Omega}$.
\end{defi}
 \begin{theorem}(\cite{Jer-Dron})(The Sobolev embedding) Let $\Omega \subset \R^d$ a weakly open Lipschitzian space.
 \begin{enumerate}
 \item[i)] If $p \in [1, d[$ then the embedding  $ W^{1-\frac{1}{p},p}(\partial\Omega)
  \hookrightarrow L^{ \frac{(d-1)p}{d-p}}(\partial\Omega)$
 is continuous.
 \item[ii)] If $p= d$ then the embedding  $ W^{1-\frac{1}{p},p}(\partial\Omega)
  \hookrightarrow L^{ q}(\partial\Omega)\, \, $ for all $ q\, \in [1, \infty[ $ is continuous.
  \item[iii)]  If $p \in ]d, \infty[$ then the embedding  $ W^{1-\frac{1}{p},p}(\partial\Omega)
  \hookrightarrow C^{ 0, 1-\frac{d}{p}}(\partial\Omega)$
 is continuous.
 \item[iv)] $ W^{1,\infty}(\partial\Omega) = C^{0,1}(\partial\Omega))$ and the norms $\|.\|_{W^{1,\infty}(\partial\Omega)}$ and $\|. \|_{C^{0,1}(\partial\Omega)}$ are equivalent.
 \end{enumerate}
 \end{theorem}
\subsection{Appendix C: Some examples of adequate open subsets}Conditions $(H1)$ contain physical reasons and technical reasons, in order to see examples, we give some in the first subsection. To simplify the idea of the technical conditions $(H1)$, we give the following examples:
\begin{enumerate}
\item[$\bullet$] {\bf Example 1:} For $R>0$ and $0<\theta_1<\theta_2<2\pi$, we consider the open domain and their partition boundary
$$\left\{\begin{array}{lll}
\Omega&=&D(0,R)\\
\Gamma_D&=&\{(R\cos(\theta),R\sin(\theta)),\;0<\theta<\theta_1\}\\
\Gamma_N&=&\{(R\cos(\theta),R\sin(\theta)),\;\theta_1<\theta<\theta_2\}\\
\Gamma_R&=&\{(R\cos(\theta),R\sin(\theta)),\;\theta_2<\theta<2\pi\}.
\end{array}\right.$$
\item[$\bullet$] {\bf Example 2:} For $R>0$, $n\in\N$, $n\geq4$ and $0<\theta_1<\theta_2<...<\theta_n=2\pi$, we consider
$A_1,\;A_2,\;A_3$ a partition of $\{1,...,n\}$:
$$\left\{\begin{array}{l}
\{1,...,n\}=A_1\cup A_2\cup A_3\\
\forall i\in\{1,2,3\}:\;A_i\neq\emptyset\\
\forall i\neq j\in\{1,2,3\}:\;A_i\cap A_j=\emptyset.
\end{array}\right.$$
The open domain and their partition boundary is defined by:
$$\left\{\begin{array}{lll}
\Omega&=&D(0,R)\\
\Gamma_D&=&\bigcup_{i\in A_1}\{(R\cos(\theta),R\sin(\theta)),\;\theta_i<\theta<\theta_{i+1}\}\\
\Gamma_N&=&\bigcup_{i\in A_2}\{(R\cos(\theta),R\sin(\theta)),\;\theta_i<\theta<\theta_{i+1}\}\\
\Gamma_R&=&\bigcup_{i\in A_3}\{(R\cos(\theta),R\sin(\theta)),\;\theta_i<\theta<\theta_{i+1}\}.
\end{array}\right.$$
\item[$\bullet$] {\bf Example 3:} Let $R>0$ and $(\theta_n)_{n\in\N}$ be a strictly monotone sequence such that:\\
If $(\theta_n)_{n\in\N}$ is a increasing sequence, we have
$$\theta_1=0\;{\rm and}\;\lim_{n\rightarrow+\infty}\theta_n=2\pi.$$
If $(\theta_n)_{n\in\N}$ is a decreasing sequence, we have
$$\theta_1=2\pi\;{\rm and}\;\lim_{n\rightarrow+\infty}\theta_n=0.$$
We consider
$A_1,\;A_2,\;A_3$ a partition of $\N$:
$$\left\{\begin{array}{l}
\N=A_1\cup A_2\cup A_3\\
\forall i\in\{1,2,3\}:\;A_i\neq\emptyset\\
\forall i\neq j\in\{1,2,3\}:\;A_i\cap A_j=\emptyset.
\end{array}\right.$$
The open domain and their partition boundary is defined by:
$$\left\{\begin{array}{lll}
\Omega&=&D(0,R)\\
\Gamma_D&=&\bigcup_{i\in A_1}\{(R\cos(\theta),R\sin(\theta)),\;\theta_i<\theta<\theta_{i+1}\}\\
\Gamma_N&=&\bigcup_{i\in A_2}\{(R\cos(\theta),R\sin(\theta)),\;\theta_i<\theta<\theta_{i+1}\}\\
\Gamma_R&=&\bigcup_{i\in A_3}\{(R\cos(\theta),R\sin(\theta)),\;\theta_i<\theta<\theta_{i+1}\}.
\end{array}\right.$$
\item[$\bullet$] {\bf Example 4:} Let $R>0$ and $(\theta_n)_{n\in\N}$ be a real sequence such that:
$$\left\{\begin{array}{l}
\inf\{\theta_n;\;n\in\N\}=0\\
\sup\{\theta_n;\;n\in\N\}=2\pi\\
{\rm the\,set\,of\,accumulation\,points\,of\,} (\theta_n)_{n\in\N}\;{\rm is\; finite}.
\end{array}\right.$$
We consider
$A_1,\;A_2,\;A_3$ a partition of $\N$:
$$\left\{\begin{array}{l}
\N=A_1\cup A_2\cup A_3\\
\forall i\in\{1,2,3\}:\;A_i\neq\emptyset\\
\forall i\neq j\in\{1,2,3\}:\;A_i\cap A_j=\emptyset.
\end{array}\right.$$
The open domain and their partition boundary is defined by:
$$\left\{\begin{array}{lll}
\Omega&=&D(0,R)\\
\Gamma_D&=&\bigcup_{i\in A_1}\{(R\cos(\theta),R\sin(\theta)),\;\theta_i<\theta<\theta_{i+1}\}\\
\Gamma_N&=&\bigcup_{i\in A_2}\{(R\cos(\theta),R\sin(\theta)),\;\theta_i<\theta<\theta_{i+1}\}\\
\Gamma_R&=&\bigcup_{i\in A_3}\{(R\cos(\theta),R\sin(\theta)),\;\theta_i<\theta<\theta_{i+1}\}.
\end{array}\right.$$
\end{enumerate}
\subsection{Appendix D: Open Problems}
We present in this section some open problems:
\subsubsection{The domain $\Omega \in \R^d$ , $d \geq 3$} ${}$\\
\noindent We consider the general Robin problem $(GRP)$ defined in the domain $\Omega  \in \R^d$ , $d \geq 3$:
$$
(GRP)\left\{\begin{array}{l}
\displaystyle\Delta u= 0\;{\rm in}\;\Omega\\
\displaystyle u= 0\;{\rm on}\;\Gamma_D\\
\displaystyle\frac{\partial u}{\partial n}= \phi(x)\;{\rm on }\;\Gamma_N\\
\displaystyle\frac{\partial u}{\partial n}+\varphi(x)u+\psi(x)u^2=g(x)\;{\rm on}\;\Gamma_R.
\end{array}\right.
$$
%-------
\subsubsection{The nonlinear Robin condition:} $\displaystyle\frac{\partial u}{\partial n}+\varphi(x)u+\psi(x)u^m =g(x) $ where  $m \geq 3 $.\\
By using the same method in the section \ref{sec-GnPrb}, we need to define a new set of admissible coefficient $\Phi_{ad,m}^{\varepsilon_1,\varepsilon_2}$  such that if $(\varphi,\psi)\in\Phi_{ad,m}^{\varepsilon_1,\varepsilon_2}$ the following problem has a unique solution in the space $V$:
$$
(RP)\left\{\begin{array}{l}
\displaystyle\Delta u= 0\;{\rm in}\;\Omega\\
\displaystyle u= 0\;{\rm on}\;\Gamma_D\\
\displaystyle\frac{\partial u}{\partial n}= \phi(x)\;{\rm on }\;\Gamma_N\\
\displaystyle\frac{\partial u}{\partial n}+\varphi(x)u+\psi(x)u^m= g(x)\;{\rm on}\;\Gamma_R.
\end{array}\right.
$$
\subsubsection{An exponential Robin condition}
We consider a more general problem where the the Robin condition is of exponential type(the practical case):
$$
(RP)\left\{\begin{array}{l}
\displaystyle\Delta u= 0\;{\rm in}\;\Omega\\
\displaystyle u= 0\;{\rm on}\;\Gamma_D\\
\displaystyle\frac{\partial u}{\partial n}= \phi(x)\;{\rm on }\;\Gamma_N\\
\displaystyle\frac{\partial u}{\partial n}+\lambda (e^{\alpha u} - e^{-(1-\alpha )u})= g(x)\;{\rm on}\;\Gamma_R.
\end{array}\right.
$$
To study this problem, we think  that is very difficult to give a set of admissible functions and the method to use is very different.

\end{document}